\documentclass{amsart} % Nicer than default article stylem: less 
\usepackage{amsmath,amsthm}     % Handy math stuff, theorem environments.
\usepackage{amssymb}            % Fancy math symbols.
\usepackage{euscript}           % Nice script font.
\usepackage{enumerate,calc}
\usepackage[matrix,arrow,curve,frame]{xy}    % XY-pic diagram pac

\xymatrixcolsep{1.9pc}                          % Adjust size of diagrams.
\xymatrixrowsep{1.9pc}
\newdir{ >}{{}*!/-5pt/\dir{>}}                  % Make better tailed arrows

\def\longbfib{\DOTSB\twoheadleftarrow\joinrel\relbar}
\def\longfib{\DOTSB\relbar\joinrel\twoheadrightarrow}

\newtheorem{thm}[subsection]{Theorem}
\newtheorem{defn}[subsection]{Definition}
\newtheorem{prop}[subsection]{Proposition}

\newtheorem{cor}[subsection]{Corollary}
\newtheorem{lemma}[subsection]{Lemma}

\theoremstyle{definition}  % Bold headings and Roman body text.

\newtheorem{remark}[subsection]{Remark}

  {\end{list}}

\newcommand{\dfn}{\textbf} % Make defined words bold.

\newcommand{\mdfn}[1]{\dfn{\mathversion{bold}#1}} % Even make math bold

\newcommand{\cat}{\EuScript}    % Use \EuScript to name a category.
\newcommand{\cC}{{\cat C}}

\newcommand{\cD}{{\cat D}}

\newcommand{\cM}{{\cat M}}

\newcommand{\Set}{{\cat Set}}

\newcommand{\sSet}{s{\cat Set}}

\DeclareMathOperator{\Hom}{Hom}

\DeclareMathOperator{\Map}{Map}

\DeclareMathOperator{\diag}{diag}

\newcommand{\ra}{\rightarrow}                   % right arrow
\newcommand{\lra}{\longrightarrow}              % long right arrow
\newcommand{\la}{\leftarrow}                    % left arrow
\newcommand{\lla}{\longleftarrow}               % long left arrow
\newcommand{\llra}[1]{\stackrel{#1}{\lra}}      % labeled long right
\newcommand{\llla}[1]{\stackrel{#1}{\lla}}      % labeled long right
\newcommand{\bwe}{\llla{\sim}}
\newcommand{\cof}{\rightarrowtail}              % cofibration
\newcommand{\fib}{\twoheadrightarrow}           % fibration
\newcommand{\trfib}{\stackrel{\sim}{\longfib}}
\newcommand{\trcof}{\stackrel{\sim}{\cof}}
\newcommand{\btrfib}{\stackrel{{\sim}}{\longbfib}}
\newcommand{\btrcof}{\stackrel{{\sim}}{\bcof}}
\newcommand{\inc}{\hookrightarrow}              % inclusion
\newcommand{\Id}{Id}                            % The identity functor
\newcommand{\id}{Id}                            % The identity functor

\newcommand{\ovcat}{\downarrow}

\numberwithin{equation}{subsection}

\newenvironment{myequation}
  {\addtocounter{subsection}{1}\begin{equation}}
  {\end{equation}$\!\!$}

\newcommand{\bcof}{\leftarrowtail}              % cofibration

\DeclareMathOperator{\WCofib}{WCofib}
\DeclareMathOperator{\WFib}{WFib}

\newcommand{\ZZ}{(\WCofib)^{-1}\cM(\WFib)^{-1}}

\DeclareMathOperator{\Ob}{Ob}

\newcommand{\rah}{\llra{h}}
\newcommand{\rav}{\llra{v}}

\begin{document}

\title{Classification spaces of maps in model categories}

%\date{February 21, 2006}
\author{Daniel Dugger}
\address{Department of Mathematics\\ University of Oregon\\ Eugene, OR
97403} 

\email{ddugger@math.uoregon.edu}

\begin{abstract}
We correct a mistake in \cite{function} and use this to identify
homotopy function complexes in a model category with the nerves of
certain categories of zig-zags.
\end{abstract}

\maketitle

%\tableofcontents

\section{Introduction}
Let $\cM$ be a model category and let $X,Y\in \cM$.  Consider the
category $\cM(X,Y)_{\Hom}$ in which an object is a zig-zag of the form
\[ X \bwe U \ra V \bwe Y, 
\]
where the indicated maps are weak equivalences.  A map from $[X\la U
\ra V \la Y]$ to $[X \la U' \ra V' \la Y]$ consists of weak
equivalences $U\ra U'$ and $V\ra V'$ making the evident diagram
commute.  We'll call $\cM(X,Y)_{\Hom}$ the \dfn{moduli category} of
maps from $X$ to $Y$.  The nerve of this category will be called the
moduli space of maps, or the {\it classification space\/} of maps.
Dwyer-Kan proved that this nerve has the correct homotopy type for a
homotopy function complex from $X$ to $Y$ 
(see \cite[6.2,8.4]{calculating}).

If $Y$ is fibrant (or, dually, if $X$ is cofibrant) one might expect 
to be able to use a simpler category here.  Namely, let
$\cM(X,Y)_{\Hom-f}$ denote the category whose objects are zig-zags
\[ X \bwe U \lra Y
\]
and where a map from $[X\bwe U \lra Y] \ra [X\bwe U' \lra Y]$ is a weak
equivalence 
$U\ra U'$ making the diagram commute.  Notice that there is an
inclusion 
\[ \cM(X,Y)_{\Hom-f}\inc \cM(X,Y)_{\Hom},
\]
as the former is
just the full subcategory consisting of objects $[X\la U \ra V \la Y]$
in which $Y\ra V$ is the identity map.

\begin{thm}
\label{th:main}
If $Y$ is fibrant then $\cM(X,Y)_{\Hom-f} \ra \cM(X,Y)_{\Hom}$ induces
a weak equivalence on nerves.
\end{thm}

Theorem~\ref{th:main} was stated in \cite[Remark 2.7]{BDG}, and was
later needed in the paper \cite{post}.  A proof was not given in
\cite{BDG}, but it is said there that the theorem follows from the
arguments in \cite[7.2]{function}.  Unfortunately, there turns out to
be a mistake in \cite[7.2]{function}, which we describe in
Remark~\ref{re:mistake} below.  The purpose of this short note is to
correct this mistake and to prove Theorem~\ref{th:main}.

We remark that the categories $\cM(X,Y)_{\Hom-f}$ seem to  be of some
current interest.  In addition to their use in \cite{BDG} and
\cite{post}, one also finds them appearing in \cite{J}.

\subsection{Notation}
Throughout the paper we will blur the distinction between a category
$\cC$ and its nerve $N\cC$.  A functor $F\colon \cC \ra \cD$ is said
to be a weak equivalence if it induces a weak equivalence on nerves.
The functor is called a homotopy equivalence if there is a functor
$G\colon \cD \ra \cC$ together with zig-zags of natural
transformations between $F\circ G$ and $\id_\cD$, and between $G\circ
F$ and $\id_\cC$.

Also, if $\cC$ is a category and $X,Y\in \Ob(\cC)$ then we write
$\cC(X,Y)$ for the set of maps from $X$ to $Y$.

\medskip

The author is grateful to Phil Hirschhorn and Dan Kan for helpful
conversations about this material.  Joachim Kock kindly provided some
corrections in terminology, as well as references for double
categories.  Also, after writing this paper the author learned from
Mandell that the Dwyer-Kan mistake has also been corrected in
\cite[Section 7]{M}.

\section{Moduli categories of maps}
Let $\cM$ be a model category and let $X,Y\in \cM$.  
Write $\WFib$ and $\WCofib$ for the categories of trivial fibrations
and trivial cofibrations.  Following \cite{function}, let
 $\ZZ(X,Y)$ denote
the full subcategory of $\cM(X,Y)_{\Hom}$ whose objects are diagrams
\[ 
\xymatrix{ X  & U \ar[r] \ar@{->>}[l]_\sim & V & Y. \ar@{ >->}[l]_\sim
}
\]

\begin{prop}
The inclusion  $\ZZ(X,Y) \inc
\cM(X,Y)_{\Hom}$ is a homotopy equivalence.
\end{prop}

\begin{proof}
Given an object $[X \la U \ra V \la Y]$
in $\cM(X,Y)_{\Hom}$, functorially factor $U\ra X$ as $U\trcof U' \fib X$.
Let $V'$ be the pushout of $U' \btrcof U \ra V$, and note that $V\ra
V'$ is a trivial cofibration.  Next, functorially factor the
composite $Y\ra V'$ as $Y\trcof V'' \fib V'$.  Let $U''$ be the
pullback of $U' \ra V' \bwe V''$.
One has the resulting
diagram
\[ \xymatrix{
X\ar@{=}[d] & U \ar[r]\ar@{ >->}[d]_\sim \ar[l]_\sim & V\ar@{ >->}[d]_\sim &
Y\ar@{=}[d] \ar[l]_\sim \\  
X & U' \ar[r] \ar@{->>}[l]_\sim & V' & Y \ar[l]_\sim \\ 
X\ar@{=}[u] & U'' \ar@{->>}[u]^\sim\ar[r] \ar@{->>}[l]_\sim &
V''\ar@{->>}[u]^\sim & Y.\ar@{=}[u] \ar@{ >->}[l]_\sim }
\]
Define a functor $F\colon \cM(X,Y)_{\Hom} \ra \ZZ(X,Y)$ by sending the object
$[X\la U \ra V \la Y]$ to $[X\la U'' \ra V'' \la Y]$.  Let $j$ denote
the inclusion in the statement of the proposition.  The above diagram
shows that there are zig-zags of natural weak equivalences between
$jF$ and the identity, and between $Fj$ and the identity.
\end{proof}

Define $\cM(\WFib)^{-1}(X,Y)$ to be the subcategory of
$\cM(X,Y)_{\Hom-f}$ whose objects are zig-zags $X \btrfib U \ra Y$.

\begin{prop} 
\label{pr:we}
The inclusion $\cM(\WFib)^{-1}(X,Y) \inc
\cM(X,Y)_{\Hom-f}$ is also a homotopy equivalence, provided $Y$ is fibrant.
\end{prop}

\begin{proof}
The proof in this case is a little different than the above.  Given a
diagram
\[ X \bwe U \lra Y\]
functorially factor the induced map $U\ra X\times Y$ as $U\trcof U'
\fib X\times Y$.   Since $Y$ is fibrant, the projection $X\times Y \ra
X$ is a fibration; so the composite $U' \ra X\times Y \ra X$ is also a
fibration (and hence a trivial fibration, by the two-out-of-three
property).  

Define $F\colon \cM(X,Y)_{\Hom-f} \ra \cM(\WFib)^{-1}(X,Y)$ by sending
$[X\la U \ra Y]$ to $[X\la U' \ra Y]$.  It is easy to check that this
gives a homotopy inverse for the inclusion.     
\end{proof}

Because of the evident commutative square
\[ \xymatrix{
\ZZ(X,Y) \ar[r]^-\sim & \cM(X,Y)_{\Hom} \\
\cM(\WFib)^{-1}(X,Y)\ar[u]\ar[r]^\sim & \cM(X,Y)_{\Hom-f}\ar[u]
}
\]
we now know that the right vertical map is a weak equivalence if and
only if the left vertical map is so.  For the rest of the paper we
will concentrate on the left vertical map.

Define another category $\cM(X,Y)_{\Hom-tw}$ in the following way.
Its objects are again 
zig-zags $X \bwe U \ra V \bwe Y$, but now a map from 
$[X \bwe U \ra V \bwe Y]$ to 
$[X \bwe U' \ra V' \bwe Y]$ is a commutative diagram of the form
\[ \xymatrix{
X \ar@{=}[d] & U \ar[r]\ar[d] \ar[l]_\sim & V & Y
\ar[l]_\sim\ar@{=}[d]\\
X  & U' \ar[r] \ar[l]_\sim & V'\ar[u] & Y.
\ar[l]_\sim
}
\]
Thus, this category has the same objects as $\cM(X,Y)_{\Hom}$ but a
different collection of morphisms.  We think of it as a `twisted'
version of $\cM(X,Y)_{\Hom}$.  

In the same way, define the category
$\ZZ_{tw}$; it is the obvious subcategory of
$\cM(X,Y)_{\Hom-tw}$.  
Note that there are inclusions 
\[ \cM(X,Y)_{\Hom-f} \inc
\cM(X,Y)_{\Hom-tw} \quad\text{and}
\]
\[
\cM(\WFib)^{-1}(X,Y) \inc \ZZ_{tw}(X,Y).
\]

Write $(\WFib\ovcat X)$ for the category of trivial fibrations with
codomain $X$.  A map in this category is a commutative triangle
\[ \xymatrix{U_1 \ar[rr]^\sim\ar@{->>}[dr]_\sim && U_2 \ar@{->>}[dl]^\sim \\
& X.}
\]
The category $(Y\ovcat\WCofib)$ of trivial cofibrations under $Y$ is
defined analogously.

\begin{remark}
\label{re:mistake}
We can now explain the mistake in \cite{function} referred to in the
introduction.  Consider the functor
\[ K\colon (\WFib\ovcat X)^{op} \times (Y\ovcat \WCofib) \ra \Set
\inc \sSet
\] 
given by $K(U\trfib X, Y\trcof V)=\cM(U,V)$.  It is
claimed in \cite[7.2iii]{function} that the homotopy colimit of this
functor (equivalently, the simplicial replacement) is isomorphic to
the nerve of $\ZZ$.  This is not correct,
however.  One readily checks that the homotopy colimit is the nerve of
$\ZZ_{tw}$ instead.
\end{remark}

Let $c\cM$ and $s\cM$ denote the categories of cosimplicial and
simplicial objects over $\cM$.  Recall that these have Reedy model
category structures, as described in \cite[Sec. 15.3]{H}.  Also,
recall that for any $Z\in \cM$ one has the associated constant
cosimplicial and simplicial objects with value $Z$; we will also
denote these $Z$, by abuse.

Let $S\colon \Delta \ra \sSet$ denote the functor $[n]\mapsto
\Delta^n$.  If $K$ is any simplicial set, let $\Delta K$ denote the
overcategory $(S\ovcat K)$---this is the {\it category of simplices\/}
of $K$.  An object of $\Delta K$ is a pair $([n],s\colon \Delta^n \ra
K)$, and the maps are the obvious ones.  We use $\Delta^{op}K$ to
denote the opposite of this category.

Note that the nerve of $\Delta K$ is homotopy equivalent to $K$.  To
see why, regard $K$ as a functor $\Delta^{op}\ra \Set \inc \sSet$.
The nerve of $\Delta K$ is the same as the simplicial replacement of
this functor, which is also the same as the homotopy colimit.  But by
\cite[Thm. 19.8.7]{H} the homotopy colimit is weakly equivalent to the
realization of this functor, which is $K$ itself.

If $QX_* \trfib X$ is a Reedy cofibrant resolution of $X$ in
$c\cM$, note that there is an evident functor
\[ \Delta \Map(QX_*,Y) \ra \cM(\WFib)^{-1}(X,Y) \]
sending $([n],QX_n \ra Y)$ to $X \btrfib QX_n \ra Y$.
Similarly, if $Y\trcof RY_*$ is a Reedy fibrant resolution of $Y$ then
there is a functor
\[ \Delta^{op} \Map(X,RY_*) \ra (\WCofib)^{-1}\cM(X,Y).
\]
The arguments of \cite{function} show the following (for the notions
of homotopy cofinal, see \cite[Def. 19.6.1]{H}):

\begin{thm} Let $QX_* \trfib X$ be a Reedy cofibrant resolution of
$X$ in $c\cM$, and let $Y\trcof RY_*$ be a Reedy fibrant resolution of
$Y$ in $s\cM$.
\begin{enumerate}[(a)]
\item The functor $Q\colon \Delta \ra (\WFib\ovcat X)$ is homotopy
left cofinal.
\item The functor $R\colon \Delta^{op} \ra (Y\ovcat\WCofib)$ is homotopy
right cofinal.
\item The map $\Delta \Map(QX_*,Y) \ra \cM(\WFib)^{-1}(X,Y)$ is a weak
equivalence. 
\item The map $\Delta^{op} \Map(X,RY_*) \ra (\WCofib)^{-1}\cM(X,Y)$ is a
weak equivalence.
\item The map $\Delta \diag\Map(QX_*,RY_*) \ra \ZZ_{tw}(X,Y)$ is a
weak equivalence.
\end{enumerate}
\end{thm}

\begin{proof}
Although this is mostly in \cite{function}, we give a brief sketch for
the reader's convenience.  

For part (a), fix a trivial fibration $U\trfib X$.  The overcategory
$(Q\ovcat [U\ra X])$ has objects pairs $([n],QX_n\ra U)$ where the
composite $QX_n\ra U\ra X$ is the same as the fixed map $QX_n\ra X$.
This is precisely the category of simplices $\Delta\Map_X(QX_*,U)$
where the maps are being computed in the overcategory $(\cM\ovcat X)$.
But by \cite[Prop. 16.4.6(1)]{H}, the map $\Map_X(QX_*,U) \ra
\Map_X(QX_*,X)$ is a trivial fibration.  Since the codomain is just a
point, we have $\Map_X(QX_*,U)$ is contractible.  But the nerve of
$\Delta\Map_X(QX_*,U)$ is homotopy equivalent to $\Map_X(QX_*,U)$, by
the remarks preceding the theorem.  Putting everything together, we
have shown that $(Q\ovcat [U\ra X])$ has contractible nerve.  The
proof for (b) is similar.

For (c), consider the functor $F\colon (\WFib\ovcat X)^{op} \ra
\Set\inc \sSet$ sending $U\ra X$ to $\cM(U,Y)$.  The homotopy colimit
(or simplicial replacement) of this functor is precisely the nerve of
the category $\cM(\WFib)^{-1}(X,Y)$.  Likewise, the homotopy colimit
of the composite
\[ \Delta^{op} \llra{Q^{op}} (\WFib\ovcat X)^{op} \llra{F} \Set \inc \sSet
\]
is the nerve of $\Delta\Map(QX_*,Y)$.  Using that the functor $Q^{op}$
is homotopy right cofinal, the induced map of homotopy colimits is a
weak equivalence by \cite[Thm. 19.6.7]{H}.  This proves (c), and (d)
is similar.

For (e), the functor in the statement is the obvious one which sends a
pair $([n],QX_n \ra RY_n)$ to the zig-zag $[X\la QX_n \ra RY_n \la Y]$
To prove that this is a weak equivalence, one introduces
 $F\colon (\WFib\ovcat X)^{op}\times
(Y\ovcat \WCofib) \ra \Set$ sending the pair $[U\ra X]$, $[Y\cof V]$
to $\cM(U,V)$.  The homotopy colimit of $F$ is the nerve of 
$\ZZ_{tw}(X,Y)$.  Now consider the composite
\[ \xymatrixcolsep{2.8pc}\xymatrix{
\Delta^{op} \ar[r]^-{\diag} & \Delta^{op}\times \Delta^{op}
\ar[r]^-{Q^{op}\times R} & (\WFib\ovcat X)^{op}
\times (Y\ovcat \WCofib) \ar[r]^-{F} & \Set.
}
\]
The homotopy colimit of the composite is the nerve of $\Delta(\diag
\Map(QX_*,RY_*))$.  The result now follows from the fact that the
functors $\diag$ and $Q^{op}\times R$ are homotopy right cofinal.
\end{proof}

\begin{cor}
\label{co:DKtw}
$\cM(\WFib)^{-1}(X,Y) \inc \ZZ_{tw}(X,Y)$ is a weak
equivalence, provided that $Y$ is fibrant.
\end{cor}

\begin{proof}
Consider the square
\[ \xymatrix{
\Delta\Map(QX_*,Y) \ar[r]\ar[d] & \cM(\WFib)^{-1}(X,Y) \ar[d] \\
\Delta\diag \Map(QX_*,RY_*) \ar[r] & \ZZ_{tw}(X,Y).
}
\]
This square does not commute, but there is a natural transformation
from one of the composites to the other---so the induced diagram of
nerves commutes in the homotopy category of simplicial sets.
  
The horizontal maps are weak equivalences by the above proposition.
The left vertical map is a weak equivalence by \cite[Prop. 17.4.6]{H},
using that $Y$ is fibrant.  So the right vertical map is also a weak
equivalence.
\end{proof}

This section has been concerned with defining terminology and
identifying exactly what is proven in \cite{function}.  The heart of
the paper lies in the next two sections, where our 
goal is to replace $\ZZ_{tw}(X,Y)$ by $\ZZ(X,Y)$ in the
above result.

%%%%%%%%%%%%%%%%%%%%%%%%%%%%%%%%%

\section{Double categories}

The notion of `double category' was introduced by Ehresmann \cite{E}.
For a topologist the name `bicategory' is much more appealing, because
their relation to categories is analagous to that of bisimplicial sets 
to simplicial sets.  In \cite{W} they are actually {\it called\/}
bicategories, but unfortunately this term has a rather different
meaning among category theorists.  Somewhat reluctantly, we will stick
to Ehresmann's original term.

\begin{defn}
A \dfn{double category} $\cC$ consists of
\begin{enumerate}[(1)]
\item A category $\cC_h$ whose maps we denote $\rah$,
\item A category $\cC_v$ with the same object set as $\cC_h$, and
whose maps we denote $\rav$,
\item A collection $S$ of squares of the form 
\[ \xymatrix{ 
\bullet \ar[r]^h \ar[d]^v & \bullet\ar[d]^v \\
\bullet \ar[r]^h  &\bullet
}
\]
which we call `bi-commutative squares'.  This collection must contain
all squares of the forms
\[
 \xymatrix{ X \ar[r]^{\id} \ar[d]_\alpha^v & X\ar[d]_\alpha^v \\
X \ar[r]^{\id}  & X
}
\qquad\text{and}\qquad
 \xymatrix{ X \ar[r]^{h}_\beta \ar[d]_{\id} & X\ar[d]_{\id} \\
X \ar[r]^{h}_\beta  & X.
}
\]
It must have the property that given two overlapping
squares as in
\[ \xymatrix{
\bullet \ar[r]^h \ar[d]^v & \bullet\ar[d]^v \ar[r]^h 
& \bullet\ar[d]^v \\
 \bullet \ar[r]^h  &
\bullet \ar[r]^h  &\bullet
}
\]
if the two smaller squares are in $S$ then so is the outer square.
Finally, it must have the analogous property in which the roles of $h$
and $v$ are switched.
\end{enumerate}
\end{defn}

\begin{remark}
The information in a category naturally fits into the first two levels
of a simplicial set, via the nerve construction.  Likewise, the
information in a double category naturally fits into the first two
levels of a bisimplicial set.  See Section~\ref{se:nerve}.
\end{remark}

Let $I$ be a small double category.  We define a functor $F\colon I\ra \cC$
to consist of two ordinary functors $F_h\colon I_h \ra \cC_h$ and
$F_v\colon I_v \ra \cC_v$ which have the same behavior on objects and
send bicommutative squares to bicommutative squares.  The collection
$\Map(I,\cC)$ of such functors itself forms a double category, in the
following way.  One defines
a morphism in $\Map(I,\cC)_h$ from $F_1$ to $F_2$ to be a collection
of $h$-maps $F_1(X) \rah F_2(X)$ such that for any $h$-map $X\ra Y$
one gets a commutative square and for any $v$-map $X\ra Y$ one gets a
bicommutative square.  Morphisms in $\Map(I,\cC)_v$ are defined
similarly, and the notion of bicommutative square is inherited in the
obvious way from $\cC$.

\subsection{Nerves of double categories}
\label{se:nerve}

Let $\Delta^n_v$ be the double category whose underlying vertical category
is
\[ 0 \rav 1  \rav \cdots \rav n
\] 
and in which all $h$-maps are the identities.  Define $\Delta^n_h$
similarly.  

Let $\cC$ be a  double category.  A \mdfn{$v$-chain of length $n$} is a
functor $\Delta^n_v \ra \cC$; that is to say, it is a
sequence of maps of the form
\[ X_0 \rav X_1  \rav \cdots \rav X_n.
\]
Write $\Map(\Delta^n_v,\cC)$ for the double category of $v$-chains of
length $n$.

Define the \dfn{nerve} of
$\cC$ be the bisimplicial set $N\cC_{*,*}$ in which the simplicial set
forming the `column' $N\cC_{*,q}$ is the usual nerve of the category
$\Map(\Delta^q_v,\cC)_h$.  It follows that the elements of
$N\cC_{p,q}$ are arrays of bicommutative squares 
\[ \xymatrix{
\bullet \ar[r]^h \ar[d]^v & \bullet \ar[r]^h\ar[d]^v & \cdots \ar[r]^h & \bullet
\ar[d]^v \\
\bullet \ar[r]^h \ar[d]^v & \bullet \ar[r]^h\ar[d]^v & \cdots \ar[r]^h & \bullet
\ar[d]^v \\
\vdots\ar[d]^v & \vdots\ar[d]^v & \vdots & \vdots\ar[d]^v \\
\bullet \ar[r]^h  & \bullet \ar[r]^h & \cdots \ar[r]^h & \bullet
}
\]
in which there are $p$ $h$-arrows in each row and $q$ $v$-arrows in
each column.

Observe that the row $N\cC_{p,*}$ is the usual nerve of the category
$\Map(\Delta^p_h,\cC)_v$.  Also, observe that the $0$th row of
$N\cC_{*,*}$ is just $N\cC_v$, and the $0$th column is $N\cC_h$.
One obtains two natural maps
\begin{myequation}
\label{eq:edges}
N\cC_h \ra \diag(N\cC) \qquad\text{and}\qquad
N\cC_v \ra \diag(N\cC).
\end{myequation}

\subsection{Trivial double categories}
\label{se:trivial}
Suppose $\cC$ is an ordinary category.  One can define a double category
$\cC_{bi}$ by setting $\cC_h=\cC_v=\cC$ and letting the bicommutative
squares be the ordinary commutative squares.  
Note that in this situation the two maps of (\ref{eq:edges}) are both
of the form $N\cC \ra \diag
(N\cC_{bi})$.  These are certainly not equal, but the following is true:

\begin{prop}
\label{pr:two=same}
In the above situation, the two maps $N\cC \ra \diag(N \cC_{bi})$
represent the same map in the homotopy category of simplicial sets.
\end{prop}

\begin{proof}
Let $f_1$ and $f_2$ be the two maps in the statement of the
proposition.  We first claim that both $f_1$ and $f_2$ are weak
equivalences.  To see this, note that the $n$th column (or the $n$th
row) of $N\cC_{bi}$ is precisely the nerve of the ordinary diagram
category $\Map(\Delta^n,\cC)$.  
However, this category is easily seen to be
homotopy equivalent to $\cC$ itself.  So every horizontal (or
vertical) map in the bisimplicial set $N\cC_{bi}$ is a weak
equivalence, and it follows from this that $f_1$ and $f_2$ are weak
equivalences.

We next claim that there is a map $\chi\colon \diag(N\cC_{bi}) \ra N\cC$ which
is a splitting for both $f_1$ and $f_2$.
To see this, note that an $n$-simplex in $\diag(N\cC_{bi})$ is a
commutative diagram
\begin{myequation}
\label{eq:array}
 \xymatrix{
X_{00} \ar[r]\ar[d] & X_{01}\ar[r]\ar[d] & \cdots \ar[r] &
X_{0n}\ar[d] \\
X_{10} \ar[r]\ar[d] & X_{11}\ar[r]\ar[d] & \cdots \ar[r] &
X_{1n}\ar[d] \\
\vdots \ar[r]\ar[d] & \vdots \ar[d]\ar[r] & \vdots \ar[r] &
\vdots \ar[d] \\
X_{n0} \ar[r] & X_{n1}\ar[r] & \cdots \ar[r] &
X_{nn}
}
\end{myequation}
(and the face operator $d_i$ just deletes the $i$th row and column
simultaneously, etc.)  Our map $\chi$ will send this $n$-simplex to
the $n$-simplex of $N\cC$ represented by
\[ X_{00} \ra X_{11} \ra X_{22} \ra \cdots \ra X_{nn}.
\]
One readily checks that this is a map of simplicial sets.
To see that it splits $f_1$ and $f_2$, just observe that $f_1$ sends a
chain 
\[ Y_0 \ra Y_1 \ra \cdots \ra Y_n \]
to the array as in (\ref{eq:array}) having this chain along each
column and all horizontal maps the identities.  The map $f_2$ is
similar, but gives an array in which all vertical maps are the
identities. 

It follows that $\chi$ is a weak equivalence, and the fact that $\chi
f_1=\chi f_2$ now shows that $f_1$ and $f_2$ represent the same map in
the homotopy category.
\end{proof}

\subsection{Reduction of double categories}
If $\cC$ is a double category then one has natural maps of simplicial sets
$N\cC_h \ra \diag N\cC$ and $N\cC_v \ra \diag N\cC$.  We are
interested in conditions forcing these maps to be weak equivalences.

For the following proposition, let $Z$ denote the double category
consisting of six objects and a single zig-zag of bicommutative
squares as indicated by
\[ \xymatrix{
0 \ar[r]^h \ar[d]^v & 1\ar[d]^v  
& 2\ar[d]^v\ar[l]_h \\
 3 \ar[r]^h  &
4   & 5 \ar[l]_h.
}
\]

\begin{prop}
\label{pr:double cat}
Let $\cC$ be a double category.  Assume that 
for each $v$-map $\alpha\colon
X\rav Y$ 
there exists an $h$-functorial zig-zag  of bicommutative squares of the form
\[ \xymatrix{
X \ar[r]^h \ar[d]^v_{\alpha} & \tilde{X}\ar[d]^v  
& Y\ar[d]^v_{\Id}\ar[l]_h \\
 Y \ar[r]^h_{\Id}  &
Y   & Y \ar[l]_h^{\Id}.
}
\]
Here $h$-functorial means that the construction is a functor 
$\Map(\Delta^1_v,\cC)_h \ra \Map(Z,\cC)_h$.
Then the evident map $N\cC_h \ra \diag(N\cC)$ is a weak equivalence of
simplicial sets.
\end{prop}

\begin{proof}
Consider the functor $F\colon \Map(\Delta^{n-1}_v,\cC)_h \ra
\Map(\Delta^n_v,\cC)_h$ sending a sequence
\[ X_0 \rav X_1 \rav \cdots \rav X_{n-1} \]
to
\[ X_0 \llra{id} X_0 \rav X_1 \rav \cdots  \rav X_{n-1}.
\]
Also consider the forgetful functor $U\colon \Map(\Delta^n_v,\cC)_h
\ra \Map(\Delta^{n-1}_v,\cC)_h$ sending
\[ X_0 \rav X_1 \rav \cdots \rav X_{n} \]
to
\[ X_1 \rav X_2 \rav \cdots \rav X_n.\]
Then $U\circ F$ is the identity, and we claim there is a zig-zag of
natural transformations between $F\circ U$ and the identity.  This
follows immediately from the hypothesis of the proposition.
So $U$ and $F$ are homotopy equivalences.

It now follows easily that every horizontal face and degeneracy in
$N\cC_{*,*}$ is a weak equivalence.  Thus, 
$N\cC_h=N\cC_{*,0} \ra \diag(N\cC)$ is a weak equivalence as well.
\end{proof}

%%%%%%%%%%%%%%%%%%%%%%%%%%%%%%%%%%%%%%%%%%%%%%%%%%%%%%%%%%%%%

\section{Application to moduli categories}

Let $\cM$ be a model category and let $X,Y\in \cM$.  Let $\cC$ be the
double category for which 
\[ \cC_h=\ZZ(X,Y) \quad\text{and}\quad \cC_v=\ZZ_{tw}(X,Y).
\]
Define the bicommutative squares to be all squares which give
commutative diagrams in $\cM$.  That is, a square involving objects
$[X\la U_i \ra V_i \la Y]$, $1\leq i\leq 4$, is called bicommutative
if it gives a commutative square when restricted to the `$U$' factors,
and also gives a commutative square when restricted to the `$V$'
factors.

\begin{lemma}
For each map $\alpha\colon
A\rav B$ in $\cC_v$ 
there exists an $h$-functorial zig-zag  of bicommutative squares of the form
\[ \xymatrix{
A \ar[r]^h \ar[d]^v_{\alpha} & \tilde{A}\ar[d]^v  
& B\ar[d]^v_{\Id}\ar[l]_h \\
 B \ar[r]^h_{\Id}  &
B   & B  \ar[l]_h^{\Id}.
}
\]
Moreover, for each map $\beta\colon X\rah Y$ in $\cC_h$
there exists a $v$-functorial zig-zag  of bicommutative squares of the form
\[ \xymatrix{
X  \ar[d]^h_{\beta} & \tilde{X} \ar[l]_v\ar[d]^h \ar[r]^v
& Y\ar[d]^h_{\Id} \\
 Y   &
Y \ar[l]_v^{\Id} \ar[r]^v_{\Id} & Y.
}
\]
\end{lemma}

\begin{proof}
Let $A=[X \la U \ra V \la Y]$ and $B=[X\la U' \ra V' \la Y]$.  
Let the components of $\alpha \colon A \rav B$ be $f\colon U\ra U'$
and $g\colon V' \ra V$.  Our zig-zag is the evident one of the form
\[ \xymatrix{
[X\la U \ra V \la Y] \ar[r]^h\ar[d]^v &
[X\la U' \ra V \la Y] \ar[d]^v &
[X\la U' \ra V' \la Y] \ar[l]_h\ar[d]^v_{\Id}\\
[X\la U' \ra V' \la Y] \ar[r]^h_{\Id}&
[X\la U' \ra V' \la Y] &
[X\la U' \ra V' \la Y]. \ar[l]_h^{\Id}
}
\]
For instance, for the middle object in the top row the map $U'\ra V$
is the composite $U' \ra V' \ra V$.  The middle vertical map is the
one with components $\Id\colon U'\ra U'$ and $g\colon V'\ra V$,  the
second map in the top row has the same two components, etc.

Similarly, let the components of $\beta\colon A\rah B$ be $p\colon U
\ra U'$ and $q\colon V\ra V'$.  Our second zig-zag is the evident one
of the form
\[ \xymatrix{
[X\la U \ra V \la Y] \ar[d]^h &
[X\la U \ra V' \la Y] \ar[l]_v\ar[r]^v\ar[d]^h &
[X\la U' \ra V' \la Y] \ar[d]^h_{\Id}\\
[X\la U' \ra V' \la Y] &
[X\la U' \ra V' \la Y] \ar[l]_v^{\Id}\ar[r]^v_{\Id} &
[X\la U' \ra V' \la Y]. 
}
\]
\end{proof}

\begin{cor}
\label{co:we}
The two maps $N\cC_h \ra \diag(N\cC_{*,*}) \la N\cC_v$ are weak
equivalences.  
\end{cor}

\begin{proof}
The fact that the first map is a weak equivalence follows from the
above lemma and Proposition~\ref{pr:double cat}.  For the second map, we
use the obvious analog of Proposition~\ref{pr:double cat} in which the
roles of $h$ and $v$ have been interchanged (and where certain
zig-zags have been replaced by zag-zigs).
\end{proof}

\begin{prop}
The map $\cM(\WFib)^{-1}(X,Y) \inc \ZZ(X,Y)$ is a weak equivalence
if $Y$ is fibrant.
\end{prop}

\begin{proof}
Recall that, in addition to the inclusion $j_1$ from the statement, one also
has the inclusion
\[  j_2\colon \cM(\WFib)^{-1}(X,Y)\inc 
\ZZ_{tw}(X,Y).
\]
One gets a resulting (non-commutative) diagram of simplicial
sets
\[ \xymatrix{
& N[\ZZ(X,Y)]\ar[dr]^\sim \\
N[\cM(\WFib)^{-1}(X,Y)] \ar[ur]^{j_1}\ar[dr]_-{j_2} && \diag(N\cC) \\
& N[\ZZ_{tw}(X,Y)]\ar[ur]_\sim 
}
\]
in which the indicated maps are weak equivalences by
Corollary~\ref{co:we}.  Moreover, one checks that if we make
$\cM(\WFib)^{-1}(X,Y)$ into a double category in the trivial way (as in
Section~\ref{se:trivial}) then we actually have a map of double categories
\[ F\colon [\cM(\WFib)^{-1}(X,Y)]_{bi} \ra \cC.
\]
The two composites in the above diagram are the same as the two
composites
\[ \xymatrix{
N[\cM(\WFib)^{-1}(X,Y)] \ar@<0.5ex>[r]^{i_1}\ar@<-0.5ex>[r]_{i_2} &
N[\cM(\WFib)^{-1}(X,Y)]_{bi} \ar[r]^-F & N\cC,
}
\]
where $i_1$ and $i_2$ are the two maps from (\ref{eq:edges}).
By
Proposition~\ref{pr:two=same}, $i_1$ and
$i_2$ represent the same map in the homotpy category of simplicial
sets, so the same is true of the two composites in our diagram.
Thus, we find that $j_1$ is a weak equivalence if and only if
$j_2$ is a weak equivalence.  But Corollary~\ref{co:DKtw} verified
that the latter is a weak equivalence.
\end{proof}

Finally, we complete the proof of our main result:

\begin{proof}[Proof of Theorem~\ref{th:main}]
This follows immediately from the above result and the diagram
immediately after
Proposition~\ref{pr:we}.
\end{proof}
%%%%%%%%%%%%%%%%%%%%%%%%%%%%%%%%%%%%%%%%%%%%%%%%%%%%%%

\bibliographystyle{amsalpha}

\end{document}